

\input amstex

\magnification 1200
\documentstyle{amsppt}

\def\vspace{\vskip.10in}
\def\C{\Bbb C}

\hsize=5.5in
\nopagenumbers
\centerline{LINEAR STRUCTURE ON CALABI-YAU MODULI SPACES}
\vskip18pt
\TagsOnRight
\centerline{Z. Ran\footnote"*"{Supported in part by NSF under DMS-9202050}}
\centerline{Department of Mathematics}
\centerline{University of California}
\centerline{Riverside, CA \ 92521 \ USA}
\centerline{ziv\@ucrmath.ucr.edu}
\vskip36pt
\baselineskip=18pt
\centerline{ABSTRACT}
\vskip12pt
{\it We show that the formal moduli space of a Calabi-Yau manifold $X^n$
carries a linear structure, as predicted by mirror symmetry.  This
linear structure is canonically associated to a splitting of the
Hodge filtration on $H^n(X)$.}

\vskip54pt
We begin by establishing terminology.  In this paper a {\it Calabi-Yau
$n$-manifold} means a compact complex manifold $X$ such that

\vspace
(i)  $X$ has trivial canonical bundle $K_X=\Omega^n_X={\Cal O}_X\;;$

\vspace
(ii)  $X$ carries no holomorphic vector fields, i.e. $H^0(\Theta_X)=0\;;$

\vspace
(iii)  $X$ is K\"ahlerian, i.e. admits some K\"ahler metric.

\vspace
(Note that (ii) is equivalent to the universal cover of $X$ having no
flat factors -- cf. [P]).  Recently such manifolds, long of interest to
mathematicians, have received a great deal of attention on the part of
physicists, because of their role in conformal field theories:  see [Y]
for a collection of papers in this vein.  One consequence of said attention
has been the emergence of a number of what might be called `physical facts'
about Calabi-Yau manifolds:  these are mathematical assertions which physicists
regard as established facts while at least some mathematicians would regard
them as likely, and often extremely intriguing, but not rigoriously proven in
the usual mathematical sense.
\vfill\eject

\pageno=2
One particularly interesting such physical fact is the assertion, common
in the physics literature (cf. [FL] and references therein), that the
moduli space of (complex structures on) Calabi-Yau manifolds, at least
in dimension $n=3$, should admit a canonical linear structure or flat
coordinates: indeed this assertion would seem to be a consequence of the
conjectured Mirror Symmetry, because the moduli of complex structures on
$X$ would correspond to the moduli of complexified K\"ahler structures
on a mirror $Y$ of $X$ if $Y$ exists, and this latter moduli is just an
open subset in a complex vector space.

\vspace
On the other hand, consideration of the case $n=2$, where $X$ is a $K3$
surface, and the usual description of the moduli of $X$ in terms of
periods, makes it seem rather unlikely that such a linear structure
can exist depending holomorphically only on the `bare' manifold $X$.

\vspace
The purpose of this paper is to prove that
a canonical
linear structure on the formal moduli on $X$ can indeed be constructed,
in agreement, apparently, with the predictions of Mirror Symmetry.  To
give a precise statement we need some more terminology.

\vspace
Let $X$ be a Calabi-Yau $n$-manifold.  We identify $H^n(X, \Bbb C)
=H^n_{DR}(X)$ as the hypercohomology of the holomorphic DeRham
complex
$$\Omega^\cdot_X:\;{\Cal O}_X{\buildrel d\over\rightarrow}\Omega^1_X
{\buildrel d\over\rightarrow}\cdots\;{\buildrel d\over\rightarrow}
\Omega^n_X\,.$$
This gives rise to the Hodge filtration $F^\cdot H^n_{DR}(X)$, which
corresponds to the stupid filtration $F^\cdot(\Omega^\cdot_X)$.  By a
{\it splitting} on $X$ we shall mean a collection of complex vector
subspaces $\{H^{p,q}\subset H^n_{DR}(X):p+q=n\}$ which split the Hodge
filtration.  We do {\it not} impose the condition that $H^{p,q}=
\overline{H^{q,p}}$.  Thus the set of all splitting on $X$ forms a
Zariski open subset of a suitable product of complex Grassmannians.
It contains a canonical element given by $\{H^{p,q}=F^p\cap
\overline{F}^q\}$, which however does not vary holomorphically with
$X$.
[In this connection, the
following problem arises naturally and is apparently unsolved.

\vspace
{\it Problem O}.  {\it Is there a natural map ${\Cal K}_{\Bbb
C}(X)\to\{splittings\}$
on the complexified K\"ahler  cone, which itself varies
holomorphically with $X$?}

\vspace
Mirror Symmetry would suggest that the answer should be affirmative:\hskip.1in
Because\break\ $\oplus H^{p,q}(X)=\oplus H^{p,p}(Y)$ for a mirror $Y$ and
the latter depends holomorphically on the complex moduli of $Y$, i.e.
the K\"ahler moduli of $X$].

\vspace
Now let $M$ be the component of $[X]$ in the moduli space.  By
a theorem of Viehweg, [V], $M$ is a quasi-projective variety.  On the
other hand let $\Cal M$ be the canonical universal formal deformation
of $X$, as constructed in [R2] (so the formal neighborhood of $[X]$
in $M$ is just $\Cal M/Aut(X)$).  By the theorem of Bogomolov-Tian-Todorov,
$\Cal M$ is smooth (this will be reproved below).  Put
$$T=T_{[X]}\Cal M=H^1(X,\Theta_X)\,,\,\Theta_X=\hbox{holomorphic tangent
sheaf.}$$

{\bf Theorem 1}.  {\it Given a split Calabi-Yau manifold $(X,\{H^{p,q}\})$,
we have a canonical isomorphism}
$$\Cal M\simeq Spf({\Bbb C}[[T^\ast]])$$

{\it Remarks}.  1.  If Problem O above were solved affirmatively, the
proof below should yield similarly a linear structure on the formal
moduli $\widetilde{\Cal M}$ of the pair ($X$, complexified complex structure
$J$).

\vspace
2.  The linear structure we obtain is a formal one and we have not
addressed questions of convergence.  Even if this were resolved and
a local linear structure on $\widetilde{\Cal M}$, say, obtained, this would
seem to depend on the particular point $(X,J)$.

\vspace
We begin the proof of Theorem 1 by recalling and amplifying some
notions and constructions from [R2] concerning higher-order deformations
and, in particular, relating them to DeRham cohomology and the Gauss-Manin
connection.

\vspace
For a topological space $X$, we denote by $X\langle m\rangle$ its $m$-fold
very symmetric product, defined as the space of all nonempty subsets of $X$
of cardinality $\le m$, with the topology induced from the Cartesian product
$X^m$ or the symmetric product $X_m$, i.e. such that we have a diagram of
topological quotients.
$$\matrix&X^m\\
\pi_m&\downarrow&\searrow\\
&X_m&\buildrel \longrightarrow \over q_m&X\langle m\rangle\,.\endmatrix$$
Given a sheaf ${\Cal L}$ of complex Lie algebras on $X$, we construct
as in [R1, R2] the associated Jacobi complex $J^\cdot_m({\Cal L})$ on
$X\langle m\rangle$.  This construction generalizes naturally to
${\Cal L}$-modules, as well as to complexes of such.  Explicitly,
let $({\Cal E}^\cdot,d^\cdot)$ be a complex of ${\Cal L}$-modules.
We define the Jacobi bicomplex $J^{\cdot\cdot}=J^{\cdot\cdot}_m({\Cal L},
{\Cal E}^\cdot)$ on $X\langle m\rangle\times X$ as follows.  First, note the
natural embedding
$$X=\Delta(X)\hookrightarrow X\langle 1\rangle\times X
\subseteq X\langle m\rangle\times X\,.$$

Set
$$J^{a,-b}=J^{-b}_m({\Cal L})\boxtimes {\Cal E}^a\quad,\; b>0\,,$$
where
$$\alignat 2
J^{-b}_m({\Cal L})&=(q_b)_\ast\big(({\Cal L}
\boxtimes\ldots\boxtimes{\Cal L})^-\big)&&\quad,\\
J^{a,-b}&={\Cal E}^a\vert_X&&\quad,\;b=0\,.\endalignat$$
Making this into a bicomplex, the horizontal differentials are given by
$$d^{a,-b}=(-1)^b\,id\otimes d^a$$
while the vertical ones are defined by the usual formula from Lie algebra
homology
$$\align \partial^{a,-b}&(t_1\!\!\times\!\cdots\!\times\! t_b{\times}e){=}\\
&{=}{(-1)^a\over{b!}}Res\Big(\!\!\sum_{\sigma\in S_b}\!sgn(\sigma)
\big(\!t_{\sigma(1)}\!\!\times\!\cdots\!\times[t_{\sigma(b-1)}\!,t_{\sigma(b)}
\!]\!\times e-t_{\sigma(1)}\!\!\times\!\cdots\!\times t_{\sigma(b-1)}
\!\!\times
\!(t_{\sigma(b)^\cdot}\!\!e)\!\big)\!\Big)\endalign$$
Note that if $\Cal E^\cdot=\C$ with trivial $\Cal L$-action, then
$$J^{\cdot\cdot}_m(\Cal L\,,\,\Cal E^\cdot)=J^\cdot_m(\Cal L)\boxtimes
\C\oplus\C_X\,.$$
(Incidentally, in this paper we will identify a bicomplex with the
associated simple complex and in particular permit a map of bicomplexes
to only preserve {\it total} degree.)

\vspace
The case we will be interested in here is where $X$ is a compact
complex manifold, $\Cal L$ is (essentially) the Lie algebra $\Theta_X$
of holomorphic vector fields, and ${\Cal E}^\cdot=\Omega^\cdot_X$ is
the holomorphic DeRham complex, on which $\Theta_X$ acts by Lie
derivative.  The resulting bicomplex $J^{\cdot\cdot}_m(\Theta_X,
\Omega^\cdot_X)$, which might be called the Jacobi-DeRham bicomplex of
$X$, looks like
$$\matrix
\Cal O_X&\rightarrow\Omega^1_X\rightarrow&\cdots&&
\qquad\cdots\Omega^n_X\\
\uparrow&&&&\uparrow\\
\Theta_X\boxtimes\Cal O_X&\rightarrow&&&\rightarrow\Theta_X
\boxtimes\Omega^n_X\\
\uparrow&&&&\uparrow\\
&\\
&\\
\uparrow&&&&\uparrow\\
(\Theta_X\boxtimes\cdots\boxtimes\Theta_X)^-\boxtimes\Cal
O_X&\rightarrow&\cdots
&\qquad\cdots&\rightarrow(\Theta_X\boxtimes\cdots\boxtimes\Theta_X)^-\boxtimes
\Omega^n_X\endmatrix$$
By the Poincar\'e lemma, we have a quasi-isomorphism
$$J^{\cdot\cdot}_m(\Theta_X,\Omega^\cdot_X)\sim J^\cdot_m(\Theta_X)
\boxtimes\C\oplus\C_X\tag1$$
In essence, this quasi-isomorphism is nothing but the Gauss-Manin
connection, as we proceed to explain.

\vspace
Define the $m$-th {\it prolongation} of the DeRham cohomology
group $H^r(X)={\Bbb H}^r(\Omega^\cdot_X)$ by the formula
(in which $\Bbb H^{\cdot,\cdot}$ denotes Kunneth components):
$$D^mH^r(X)=\Bbb H^{0,r}\big(X\langle m+1\rangle\,,\,J^{\cdot\cdot}_m
(\Theta_X\,,\,\Omega^\cdot_X)\big)\,.$$
This terminology is justified by the following result, which
can be proven by the method of [R2] (proof omitted).

\vspace
{\bf Theorem 2}.  {\it Let $X_m/R_m$ be the canonical $m$-th order
deformation of $X$ as in} [R2].  {\it Then we have a canonical
isomorphism}
$$D^mH^r(X)\simeq {Diff}^m_{R_m}\big(H^r_{DR}(X_m/R_m)^\ast\,,
\C\big)\,.$$

Note that these groups form a tower
$$H^r(X)=D^0H^r(X)\subseteq D^1H^r(X)\subseteq\cdots\subseteq D^m H^r
(X)$$
and we have exact sequences
$$0\to D^{m-1}H^r(X)\rightarrow D^mH^r(X)\rightarrow S_mT\otimes H^r
(X)$$

By (1), we have
$$D^mH^r(X)\simeq T^{(m)}R_m\otimes H^r(X)\Theta H^r(X)={Diff}^m(R_m,\C)\otimes
H^r(X)\,,$$
where $T^{(m)}R_m={Diff}^m_+(R_m,\C)=\Bbb H^0\big(X\langle m\rangle,
J^\cdot_m(\Theta_X)\big)$, reflecting the Gauss-Manin connection
on $H^r_{DR}(X_m/R_m)$.  This isomorphism can also be seen, even
on the homotopy level, directly from the bicomplex $J^{\cdot\cdot}_m
(\Theta_X,\Omega^\cdot_X)$, using the Cartan formula the Lie
derivative of differential forms:
$$L_t=i_t\circ d+d\circ i_t$$
where $t\in\Theta_X$ and $i_t$ denotes interior multiplication by $t$.
Indeed we may define a homotopy splitting $\gamma$ of the natural
projection
$$J^{\cdot\cdot}_m(\Theta_X,\Omega^\cdot_X)\rightarrow J^{\cdot\cdot}
_m(\Theta_X,\Omega^\cdot_X)/G^0J^{\cdot\cdot}_m(\Theta_X,\Omega^\cdot
_X)\,,$$
where $G^\cdot$ denotes the vertically stupid filtration on a bicomplex and
the RHS is viewed as all the terms of vertical degree $<0$ in $J^
{\cdot\cdot}_m(\Theta_X,\Omega^\cdot_X)$, by the formula
$$\align \gamma^{a,-b}&=id\qquad\qquad b\ge2\\
a^{a,-1}&=id\oplus j\endalign$$
where $j:\Theta_X\boxtimes\Omega^a_X\longrightarrow\Omega^{a-1}_X$ is
restriction followed by interior multiplication.

\vspace
The cohomology map associated to $\gamma$ yields the Gauss-Manin
splitting of the inclusion
$$H^r(X)\subseteq D^mH^r(X)\,.$$
A similar construction yields a homotopy equivalence
$$J^{\cdot\cdot}_m(\Theta_X,\Omega^\cdot_X)/G^0J^{\cdot\cdot}_m(\Theta_X,
\Omega^\cdot_X)\sim J^\cdot_m(\Theta_X)\boxtimes\Omega^\cdot_X$$
so that $D^mH^r(X)\simeq H^r(X)\oplus(T^{(m)}R_m)\otimes H^r(X)$.

\vspace
Now let $X$ be a Calabi-Yau manifold, and fix some holomorphic
volume form $\Phi$.  Let $\widehat{\Theta}_X\subset \Theta_X$ denote
the subsheaf of {\it divergence-free} vector fields, i.e. those
annihilating $\Phi$ via Lie derivative.  As $\Phi$ is unique
up to a constant $\widehat{\Theta}$ is independent of the choice
of $\Phi$.  As in [R2], there is a canonical formal
moduli $\widehat{\Cal M}$ for the pair $(X,\Phi)$ and we have
$$T^{(m)}\widehat{\Cal M}=\Bbb H^0\big(J^\cdot_m(\widehat{\Theta})\big)\,.$$

Next we replace the DeRham complex $\Omega^\cdot_X$ by its quasi-isomorphic
subcomplex $\Omega^\cdot_{X,0}$ defined by
$$\alignat 2
\Omega^i_{X,0}&=\Omega^0_X&&\quad\;i\le n-2\\
&=\widehat{\Omega}^{n-1}_X&&,\quad i=n-1\qquad\hbox{(i.e. the {\it closed} }
(n-1)\hbox{-forms)}\\
&=0&&,\quad i=n\,.
\endalignat$$
This forms a complex of $\widehat{\Theta}$-modules and as above one may form
the Jacobi bicomplex $J^{\cdot\cdot}_m(\widehat{\Theta}_X,\Omega^\cdot_{X,0})$
and establish a canonical isomorphism
$$\Bbb H^r(X\langle m+1\rangle,J^{\cdot\cdot}_m(\widehat{\Theta}_X,\Omega^\cdot
_{X,0})=\widehat{D}^m H^r(X)={Diff}^m_{\widehat{R}_m}\big(H^r_{DR}(\widehat{X}
_m/\widehat{R}_m)^\ast,\Bbb C\big)$$
where $\widehat{X}_m=X_mx_{R_m}\widehat{R}_m$
is the canonical $m$-th order
deformation of $(X,\Phi)$.  Note that by Hodge theory (resp. smoothness
of $\widehat{\Cal M}$), both hypercohomology spectral sequences for the
bicomplex
$J^{\cdot\cdot}_m(\widehat{\Theta}_X,\Omega^\cdot_{X,0})$ degenerate at
$E_1$.

\vspace
Now note that interior multiplication by $\Phi$ yields an isomorphism
$$\widehat{\Theta}_X\simeq\widehat{\Omega}^{n-1}_X\,,$$
and by the Cartan formula we have
$$i_{[t_1,t_2]}\Phi=d(i_{t_1\wedge t_2}\Phi)\qquad t_1,t_2\in
\widehat{\Theta}\,.$$
This implies that the entire complex $J^\cdot_{m+1}(\widehat{\Theta})[-n+1,
-1]$, viewed vertically and pulled back to $X\langle m\rangle\times X$,
is isomorphic to a direct summand, viz. the `fully alternating' part of
the subcomplex
$$F^{n-1}J^{\cdot\cdot}_m(\widehat{\Theta}_X,\Omega^\cdot_{X,0})\,,\qquad
F^\cdot=\hbox{horizontally stupid filtration.}$$

Now to complete the proof of Theorem 1 it will suffice, by the
results of [R2], to construct canonical--in terms of the given data $(X,
\{H^{p,q}\})$--isomorphisms
$$T^{(m)}\Cal M\simeq\overset m\to{\underset {i=1}\to\oplus}S_iT\,.$$
As
$$\widehat{T}=H^1(\widehat{\Theta}_X)=\Phi^{-1}F^{n-1}H^n(X)=\Phi^{-1}
(H^{n,0}\oplus H^{n-1,1})=\Phi^{-1}H^{n,0}\oplus T\,,$$
it will suffice to construct suitable isomorphisms
\vspace
{\parindent.3in\item{$(*)_m$}}\hfill $\displaystyle{\widehat{T}^{(m)}:=T^{(m)}
\widehat{\Cal M}\simeq\overset m\to{\underset
{i=1}\to\oplus}S_i\widehat{T}\,.}$\hfill\break\

\noindent We will construct these by induction on $m$, together with
isomorphisms
\vspace
{\parindent.3in\item{$(\ast\ast)_{m-2}$}}\hfill $\displaystyle
{\widehat{D}^{m-2}H^n(X)\simeq\overset m-2\to{\underset
{i=0}\to\oplus}\overset m\to{\underset{p=0}\to\oplus}S_i\widehat{T}\otimes
H^{p,q}\,.}$\hfill\break\

\noindent For $m=1$ there is nothing to prove.  For $m=2,(\ast\ast)_0$
is already
given so it suffices to construct $(\ast)_2$.  To this end, define a map of
bicomplexes
$$\varphi^\cdot_2:J^\cdot_2(\widehat{\Theta})[-n+1,-1]\to
J^{\cdot\cdot}_0(\widehat{\Theta},\Omega^\cdot_{X,0})=\Omega^\cdot_
{X,0}$$
by
$$\align
\varphi^{-1}_2:J^{-1}_2(\widehat{\Theta})&\to\Omega^{n-1}_
{X,0}=\widehat{\Omega}^{n-1}_X\\
\varphi^{-1}_2(t)&=i_t(\Phi)\,,\\
\varphi^{-2}_2:J^{-2}_2(\widehat{\Theta})&\to
\Omega^{n-2}_{X,0}=\Omega^{n-2}_X\\
\varphi^{-2}_2(t_1\times t_2)&=i_{t_1\wedge t_2}(\Phi)\,.
\endalign$$
On cohomology, this yields a diagram
$$\matrix
0&\to&\widehat{T}&\to&\widehat{T}^{(2)}&\to&S_2\widehat{T}&\to&0\\
&\\
&&\wr\vert&&\downarrow&&\downarrow\\
&\\0&\to&F^{n-1}H^n(X)&\to&H^n(X)&\to&H^n(X)/F^{n-1}H^n(X)&\to&0\,.
\endmatrix
$$
As we are given a splitting of the bottom row and the left vertical arrow
is an isomorphism, we get a splitting of the top row.  Note that for
$n=3$ the right vertical arrow is the so-called Yukawa coupling.

\vspace
Now inductively, assuming $(\ast)_m$ and $(\ast\ast)_{m-2}$ are done, we
firstly obtain $(\ast\ast)_{m-1}$ (and even $(\ast\ast)_m$) from the
Gauss-Manin isomorphism
$$\widehat{D}^iH^n(X)\simeq H^n(X)\oplus\widehat{T}^i\otimes H^n(X)\,.$$
To obtain $(\ast)_{m+1}$, define a morphism of bicomplexes
$$\varphi^\cdot_{m+1}:J^\cdot_{m+1}(\widehat{\Theta})[-n+1,-1]\to
J^{\cdot\cdot}_{m-1}(\widehat{\Theta},\Omega^\cdot_{X,0})$$
by
$$\align
\varphi^{-j}_{m+1}:J^{-j}_{m+1}(\widehat{\Theta})&\to J^{n-1,-j+1}
_{m-1}(\widehat{\Theta},\Omega^\cdot_{X,0})\qquad i\le m\\
\varphi^{-j}_{m+1}(t_1\times\cdots\times t_j)&=\sum^j_{k=1}
(-1)^kt_1\times\cdots\times\widehat{t}_k\times\cdots\times
t_j\times i_{t_j}(\Phi)\,,\\
\varphi^{-m-1}_{m+1}:J^{-m-1}_{m+1}(\widehat{\Theta})&\to
J^{n-2,-m+1}_{m-1}(\widehat{\Theta},\Omega^\cdot_{X,0})\\
\varphi^{-m-1}_{m+1}(t_1\times\cdots\times t_{m+1})&=\sum_{i,j}(-1)
^{i+j}t_1\times\cdots\times\widehat{t}_j\times\cdots\times \widehat{t}_j
\times\cdots\times t_{m+1}\times
i_{t_i\wedge t_j}(\Phi)\,.\endalign$$
Taking cohomology, we get a diagram
$$\matrix
0&\to&\widehat{T}^{(m)}&\to&\widehat{T}^{m+1}&\to&S_{m+1}\widehat{T}&
\to&0\\
&\\
&&\downarrow&&\downarrow&&\downarrow\\
&\\
0&\to&\widehat{D}^{m{-}1}\!F^{n{-}1}\!H^n(X)&\to&\widehat{D}^{m{-}1}\!H^n\!(X)
&\to&\widehat{D}^{m{-}1}\!\big(\!H^n\!(X)\!/\!F^{n{-}1}\!H^n\!(X)\!\big)
&\to&0\,.
\endmatrix$$
By induction, we have a splitting of the bottom row and, what's more, of the
left vertical arrow.  This yields a splitting of the top row, completing the
induction step.

\vspace
{\it Remarks. 1.}  A similar construction can be used to give a simple a priori
proof of the smoothness of $\widehat{\Cal M}$ or $\Cal M$,
i.e. the degeneration at
$E_1$ of the hypercohomology spectral sequence of $J^\cdot_m(\widehat
{\Theta}_X)$: the point is that the coboundary maps factor through the exterior
derivative, which induces the zero map on cohomology.  (That the
degeneration is equivalent to unobstructedness of the corresponding
moduli problem is a general fact, implicitly proven in [R2].)

{\it 2.} For $X$ symplectic, i.e. an $Sp(n)$-manifold, it is interesting to
compare
the above construction with that of 'formal twistors' in [R1]: the latter
lifts an {\it isotropic} vector $\alpha \in H^1(X,\Theta)$  to a formal
deformation,
and, unlike the former, is manifestly (canonical and) holomorphic in
$(X,\alpha)$.
However, at least in the K3 case, the two do, in fact agree, which suggests
that
they should agree in general (in the (symplectic,isotropic) case). We hope to
return
to this point elsewhere.
\vspace
{\it Note:}  After this work was done, we were informed by Professor M.
Green that he and Professor P. Griffiths had also obtained results on
this problem.  We are grateful to Professor Green for his subsequent
incisive comments on the manuscript.

\vfill\supereject

\magnification 1200
\pageno=10

\hsize=5.5in
\centerline{REFERENCES}

\vskip.5in
\item{[FL]}  S. Ferrara and J. Louis:  ``Picard-Fuchs equations
and flat holomorphic connections from $N=2$ supergravity", in [Y],
pp. 301--315.

\vspace
\item{[P]}  S\'eminaire Palaiseau:  ``G\'eom\'etrie des surfaces
$K3$:  modules et p\'eriod\'es", {\it Ast\'erisque} {\bf 126}
(1985).

\vspace
\item{[R1]}  Z. Ran:  ``Derivatives of moduli", {\it Internat. Math.
Research Notices} (1993), 93--106.

\vspace
\item{[R2]}  \underbar{\hskip.25in}:  ``Canonical infinitesimal
deformations", preprint.

\vspace
\item{[V]}  E. Viehweg:  ``Weak positivity and stability of certain
Hilbert points", I, II, III, {\it Invent. Math.}.

\vspace
\item{[Y]}  S.T. Yau, ed:  ``Essays on mirror manifolds", {\it International
Press}, Hong Kong, (1992).

\vfill\supereject
\bye